\documentclass[12pt]{amsart}

\usepackage{amsmath, amssymb}
\newtheorem{theorem}{Theorem}
\newtheorem{conjecture}[theorem]{Conjecture}
\newtheorem{corollary}[theorem]{Corollary}
\newtheorem{proposition}[theorem]{Proposition}
\newtheorem{lemma}[theorem]{Lemma}

\newtheorem{claim}{Claim}

\def\D{\mathcal{D}}

\def\bth{\begin{theorem}}
\def\eth{\end{theorem}}

\def\bc{\begin{corollary}}
\def\ec{\end{corollary}}

\def\bcj{\begin{conjecture}}
\def\ecj{\end{conjecture}}

\newcommand{\dist}{{\rm dist}}
\newcommand{\bizveg}{{\hfill $\Box$}}

\title[The diameter game]{The diameter game}

\author[Balogh]{J\'ozsef Balogh}\address{Department of Mathematical
Sciences, University of Illinois}
\thanks{The first author's research is partially supported by NSF
grants DMS-0302804,  DMS-0603769
 and DMS-0600303, UIUC Campus Research Board grants \#06139 and p \#07048, and by OTKA grants T034475 and T049398.}
\email{jobal@math.uiuc.edu}

\author[Martin]{Ryan Martin}\address{Department of Mathematics,
Iowa State University}
\thanks{The second author's research is partially supported by
NSA grant H98230-05-1-0257.}\email{rymartin@iastate.edu}

\author[Pluh\'ar]{Andr\'as Pluh\'ar}\address{University of Szeged,
Department of Computer Science}
\thanks{The third author's research is partially supported by OTKA
grants T034475 and T049398.}\email{pluhar@inf.u-szeged.hu}



\begin{document}

\begin{abstract}
A large class of  Positional Games are defined on the complete graph
on $n$ vertices. The players, Maker and Breaker, take the edges of
the graph in turns, and Maker wins iff his subgraph has a given --
usually monotone -- property.  Here we introduce the $d$-diameter
game, which means that Maker wins iff the diameter of his subgraph
is at most $d$. We investigate the biased version of the game; i.e.,
when the players may take more than one, and not necessarily the
same number of edges, in a turn. Our main result is that we proved
that the $2$-diameter game has the following surprising  property:
Breaker wins the game in which each player chooses one edge per
turn, but Maker wins as long as he is permitted to choose $2$ edges
in each turn whereas Breaker can choose as many as
$(1/9)n^{1/8}/(\ln n)^{3/8}$.

In addition, we investigate $d$-diameter games for $d\ge 3$. The
diameter games are strongly related to the degree games. Thus, we
also provide a generalization of the fair degree game for the biased
case.

\end{abstract}

\maketitle

\section{Introduction}

In the most general setting, a positional game can be viewed as a game
in which two players -- Maker and Breaker -- occupy vertices in a
hypergraph.  Maker wins if he can occupy all vertices of one
hyperedge, otherwise    Breaker wins, i.e.,  if Breaker can occupy at
least one vertex of each hyperedge. For more information on these
games, see the excellent survey included in \cite{beckcpc05}.

\subsection{Biased positional games}

Formally a Maker-Breaker Positional Game is defined as follows.
Given an arbitrary hypergraph
${\mathcal H}=(V({\mathcal H}), E({\mathcal H}))$, Maker and Breaker
take $a$ and $b$ elements of $V({\mathcal H})$ per turn. Maker wins
by taking all elements of at least one edge $A \in E({\mathcal H})$,
otherwise Breaker wins.  In this
paper we let Maker start the game. We call such games $(a:b)$-games.
If $a=b$, the game is \textit{fair}, otherwise it is
\textit{biased}. If $a=b>1$, the game is \textit{accelerated}.

In general it is hard\footnote{In fact it is PSPACE-complete, see
\cite{Schaefer}.} to decide who wins a certain Maker-Breaker game.
The first general results were obtained by Hales and Jewett, who
formalized the strategy stealing idea of Nash, and gave a
K\"onig-Hall type criterion for Breaker's win, see \cite{H-J}. The
celebrated Erd\H{o}s-Selfridge Theorem gives a condition on the
hypergraph which guarantees the existence of a winning strategy for
the Breaker when $a=b=1$, see \cite{ES}. We shall use its general
version that is due to J\'ozsef Beck, see \cite{Beck1}.

\begin{theorem}[Beck~\cite{Beck1}] \label{ESB}
Let $E(\mathcal H)$ be   the family of winning sets of a positional
game.  Breaker has a winning strategy in the $(a:b)$-game if
\begin{equation} \sum_{A\in {E(\mathcal H)}}(1+b)^{1-|A|/a}< 1 .
\label{thm:es}
\end{equation}
\end{theorem}

\noindent\textbf{Remark.} In using Theorem~\ref{thm:es}, we say the winning strategy guaranteed by the theorem is \textbf{playing as ESB-Breaker} (after Erd\H{o}s-Selfridge and Beck).  For several games we describe, the player we denote as Maker will play as ESB-Breaker.

Note that Theorem~\ref{ESB} is best possible in the sense that for
all $a, b \in {\mathbb N}$ there is a hypergraph $\mathcal H$ with
equality in \eqref{thm:es}, and Maker having a winning strategy. In
order to appreciate the significance of biased games let us discuss
some games defined on graphs.

\subsection{Graph games}

In Shannon's switching game, Maker and Breaker take the edges of
a given graph $G$, and the winning sets are the spanning trees. The
classical theorem of Lehman  \cite{Lehman} states that in a $(1:1)$-game Maker (as a second player) wins iff $G$ contains two edge-disjoint
spanning trees.

Let the edges of $K_n$ (the complete graph on $n$ vertices) be the
``board.''  By Lehman's Theorem Maker wins   the Shannon switching
game for $n\ge 4$, so Chv\'atal and Erd\H{o}s asked \cite{Chvatal}
what is the outcome of the $(1:b)$-game? More precisely, they looked
for the {\em breaking point}, that is the smallest value of
$b_0=b_0(n)$ such that Breaker wins the $(1:b_0)$-game. It turned
out that $b_0 =\Theta (n/\ln n)$.

In general, one fixes a (usually monotone) graph property $\mathcal P$,
and Maker wins if $\mathcal P$ holds for the subgraph of his edges. Again,
the question is: What is the smallest $b_0$, for which Breaker wins the
$(1: b_0)$-game? Along this line a number of beautiful results are published: for Hamiltonicity and maximum degree see Beck, \cite{BeckHam, Beck2}, for
planarity, colorability and graph minor games, see Hefetz et al,
\cite{Hefetz}, for building a specific graph $G$ or creating a large component, see Bednarska and {\L}uczak, \cite{luczak0, luczak}.

Here, we work with the monotone property ${\mathcal P}_d$, that is
the graph has diameter at most $d$. We denote the corresponding
$d$-diameter game by $\D_d(a:b)$, or more briefly, by $\D_d$ if
$a=b=1$.

First, we investigate the $2$-diameter game and in particular we
observe that in the $(1:1)$-game, Breaker wins.

\begin{proposition} \label{2-(1:1)}
For $n\leq 3$  Maker has and
 for $n\geq 4$  Breaker has a winning strategy in the game
$\D_2$.
\end{proposition}

Our main result is observing (and proving) the following phenomenon:
 When Maker
takes two edges in each turn, then,
 rather surprisingly, Breaker loses the game not only if his bias
 is  $2$ but even if it is $\frac{1}{9}n^{1/8}/(\log n)^{3/8}$, i.e., a little acceleration of the game  drastically changes the breaking point of Breaker. Our result is the first where such an unusual circumstance is proved.

\begin{theorem} \label{2-(2:b)}
Maker wins the game $\D_2(2:\frac{1}{9}n^{1/8}/(\log n)^{3/8})$, and
Breaker wins the game $\D_2(2:(2+\epsilon)\sqrt{n/\ln n})$ for any
$\epsilon >0$, provided $n$ is large enough.
\end{theorem}



We prove similar (but weaker) results for the game $D_d$ for
$d\ge 3$.

\begin{theorem} \label{d-(1:b)}
For any fixed $d\geq 3$, and $n$ large enough, Maker wins the
game $\D_d(1: (2d)^{-1} (n/ \ln n)^{1-1/\lceil d/2\rceil})$.

Furthermore, for every integer $a>1$ and integer $d\geq 3$, there exist $c_2=c_2(d)>0$
and $c_3=c_3(a,d)>0$ such that Breaker wins the games $\D_d(1:c_2n^{1-1/(d-1)})$ and
$\D_d(a:c_3n^{1-1/d})$, provided $n$ is big enough.
\end{theorem}

The careful reader can verify that the computations that confirm the results of Theorem~\ref{d-(1:b)} will hold even if $d$ is going to infinity with $n$, provided  $d\leq c_1\ln
n/\ln\ln n$ for some positive constant $c_1$.

 We conjecture that in the games $\D_d$ (for $d\ge 2$) the
correct breaking point is close (up to a polylog factor) to the
``Breaker's'' bound.



\subsection{Probabilistic intuition}

An important guide to understand such games is the so-called {\em
probabilistic intuition}, for more details and examples see
\cite{Beck2}. In the probabilistic intuition, we substitute the
perfect players with ``random players." In the simplest example,
when $a=b=1$, this process leads to a random two coloring of a
hypergraph by flipping a fair coin. The condition of
Theorem~\ref{ESB} says that the expected number of monochromatic
sets is smaller than one. That is the {\em first moment method}
tells the random Breaker has a chance to win, and Theorem~\ref{ESB}
says that deterministic Breaker also wins.


Indeed, the theories of {\em probabilistic methods} and positional
games are in close relation. One can find the game theoretic analogy
of the {\em second moment method} or {\em Lov\'asz' Local Lemma},
see some examples in \cite{Beck1, Beck2, PA1}. One of the most
far-reaching ideas is to model a biased game by a random graph. It
means that we consider the random graph $G(n,p)$ with $p=a/(a+b)$;
i.e., the expected number of edges in the random graph is the same as
the number of Maker's edges in the game at its conclusion. It was demonstrated
in a number of cases that Maker wins
if $\mathcal P$ holds, while Breaker wins if ${\mathcal P}$ does not
hold for $G(n,p)$ almost surely, see  \cite{Beck2, Hefetz, luczak0,
luczak}.

Note that acceleration of a game can have profound effects. Based on
the previous heuristic, one might tempted to think that the outcome
of a $(1:1)$-game is more or less the same as that of an $(a:a)$-game on the same $\mathcal H$ for $a>1$ if the winning sets are
big\footnote{Of course if $|A| \leq a$ for an $A \in E(\mathcal H)$,
then Maker wins instantly.}. A counterexample to this belief is the
$k$-in-a-row game. The $(1:1)$-game is a draw (i.e., Breaker wins)
for  $k \geq 8$, see \cite{Zetters}. The $(n:n)$-game is a Maker's
win not only for, say, $k \leq n +8$ but as long as $k \leq n+\log
n/\log \log n$ and it is a Breaker's win for $k \geq n +c \log n$,
see \cite{PA2}.

Since the diameter of $G(n,1/2)$ is $2$ almost surely,
we can interpret Proposition~\ref{2-(1:1)} that the
probabilistic intuition fails completely in the $\D_2(1:1)$-game.
Furthermore Theorem~\ref{2-(2:b)} means that the probabilistic
intuition is at least partially restored in the $(2:2)$-game.
The change is more dramatic than in the case of $k$-in-a-row game,
since once Maker can mark two edges per turn, the breaking point
jumps up to at least $\frac{1}{9}n^{1/8}/(\log n)^{3/8}$. Probably,
$b_0$ is much larger, perhaps even $\Omega(\sqrt{n/ \log n})$.

It is worth discussing the problem $\D_3(1:b)$ before the general
case. Observe that the threshold of the property
${\mathcal P}_3$ is about $n^{-2/3}$; i.e., $G(n,p)$ has property
${\mathcal P}_3$ with probability close to $1$ if
$p =n^{-2/3+\epsilon}$, and it does not have property ${\mathcal P}_3$
if $p=n^{-2/3-\epsilon}$ for arbitrary $\epsilon>0$ and $n$ sufficiently large.
See Bollob\'as~\cite{BBrg} for a detailed account of the diameter of
the random graph. According to Theorem~\ref{d-(1:b)},
Maker wins the game $\D_3(1:c_1\sqrt{n/\ln n})$, and Breaker wins
the game $\D_3(1:c_2\sqrt{n})$, for some $c_1, c_2>0$, provided
$n$ is big enough. That is, the game
$\D_3(1:b)$ defies the probabilistic intuition again. However we
suspect that the $\D_3(3:b)$-game again agrees with the
probabilistic intuition; i.e., the breaking point should be
$b_0\approx n^{2/3}\times {\rm polylog}(n)$. We do not have the
courage to state a conjecture for the breaking point for the
$\D_3(2:b)$-game.

In  general, the threshold function of the property ${\mathcal P}_d$
is about $n^{-(d-1)/d}$. This suggests that in the game $D_d$  the
breaking point should be around $n^{1-1/d}\times {\rm polylog}(n)$,
however for $a=1$ the probabilistic intuition still fails, see
Theorem \ref{d-(1:b)}.

The rest of the paper is organized as follows: In Section~2 we
describe general theorems and auxiliary games, such as the degree
game and expansion game. In Section~3 we prove  our results on the
$2$-diameter game and Section~4  on the $d$-diameter game for $d
\geq 3$.

\section{Auxiliary games}


We shall need a theorem of Erd\H{o}s and Chv\'atal that they used to
derive bounds for the biased Shannon switching game in \cite{Chvatal}.

\begin{theorem}[Chv\'{a}tal-Erd\H{o}s~\cite{Chvatal}]\label{thmchvatal}
Let  $\mathcal H$ be an $r$-uniform family of $k$ {\bf disjoint}
winning sets. Then
 $\!$ \\ (i)  Maker has a winning
strategy in the $(a:1)$-game when
$$ r\leq (a-1)\sum_{i=1}^{k-1}\frac{1}{i} . $$
(ii)  Maker has a winning strategy in the $(a:2)$-game when
$$ r\leq \frac{a-1}{2}\sum_{i=1}^{k-1}\frac{1}{i}.$$
\end{theorem}
In the rest of this section we state and prove theorems regarding minimum-degree and expansion games.

\subsection{Biased degree games}

In our proofs  we use  some auxiliary degree games.  Given a graph
$G$ and a prescribed degree $d$, Maker and Breaker play an $(a:b)$-game on the edges of $G$.  Maker wins by getting at least $d$ edges
incident to each vertex. For $G=K_n$ and $a=b=1$ this game was
investigated thoroughly in \cite{Szekely} and \cite{Beck2}. It was
shown that Maker wins if $d<n/2-\sqrt{n\log n}$, and Breaker wins if
$d>n/2-\sqrt{n}/12$.

This is in agreement with the probabilistic intuition, since in
$G(n,1/2)$ the degrees of all vertices fall into the interval
$\left[n/2-\sqrt{n \log n},n/2+\sqrt{n \log n}\right]$ almost surely.  When $a\neq b$, playing on $G=K_n$ analogously one would expect that
Maker wins if $d < an/(a+b)-c'\sqrt{n\log n}$, and Breaker wins if
$d > an/(a+b)-c''\sqrt{n}$ for some $c', c''>0$.

Here we are interested in giving conditions for Maker's win only,
so this will suffice:

\begin{lemma} \label{bdegree}
Let $a\leq n/(4\ln n)$ and $n$ be large enough. Then Maker wins
the $(a:b)$-degree game on $K_n$ if
$d<\frac{a}{a+b}n-\frac{6ab}{(a+b)^{3/2}}\sqrt{n \ln n}$.
\end{lemma}

\noindent\textbf{Remark.} In referencing this game, we call it the $(a:b)$-${\rm MINDEG}(d)$ game and, where the parameter is understood, refer to the winning strategy as \textbf{playing as MINDEG-Maker}.

Note that, in the statement of the lemma, no bound on $b$ is necessary.  It can be shown that if $b\geq n/(36\ln n)+a$, then the lemma simply states that Maker can achieve minimum degree at least $0$, which is trivially true. \\

\noindent {\bf Proof of Lemma~\ref{bdegree}.} We use a little
modification of the weight function argument of Beck, see
\cite{Beck2}. Consider the hypergraph ${\mathcal H}=(V({\mathcal
H}), E({\mathcal H}))$, where $V({\mathcal H})=E(K_n)$ and
$E({\mathcal H})$ contains the set $A_v$ for each vertex $v \in
K_n$, where $A_v$ is the set of edges incident to $v$.

Let $X_i$ and $Y_i$ be the set of (graph)
edges selected by Maker and Breaker, respectively, before Maker makes his last, $a^{\rm th}$ move in the $i^{\rm th}$ step. The $i^{\rm th}$ step is formally defined to be Maker's $i^{\rm th}$ set of moves, together with the succeeding moves of Breaker. The
\textit{$(\lambda_1,\lambda_2)$-weight} of a hyperedge $A$ at the end of the $i^{\rm
th}$ step is defined as
$$ w_i(A)=(1+\lambda_1)^{|Y_{i} \cap A|-(bn/(a+b)+k)}(1 -
\lambda_2)^{|X_i \cap A| - (an/(a+b) - k)}, $$ where we let
$k=\frac{6ab}{(a+b)^{3/2}}\sqrt{n \ln n}$ and the values of
$\lambda_1$, $\lambda_2$ will be given later. For any (graph)
edge $e$, let
$$ w_i(e)=\sum_{e\in A}w_i(A) \qquad\mbox{and} \qquad
T_i=\sum_{A\in E({\mathcal H})} w_i(A).$$ We want to ensure the
following three properties:
\renewcommand{\theenumi}{\roman{enumi}}
\begin{enumerate}
\item If Breaker wins the game in the $i^{\rm th}$ step, then $T_i\geq 1$, \label{it:prop1}
\item $T_{i+1}\leq T_i$, \label{it:prop2}
\item $T_0<1$. \label{it:prop3}
\end{enumerate}

Property (\ref{it:prop1}) is trivially true if Maker follows the
greedy strategy, that is he always chooses the maximum weight edge
available. Let $w$ be the weight of the largest weighted edge
\textbf{before} Maker makes his last, $a^{\rm th}$, move. This means
that Maker will reduce the value of $T_i$ by at least $a\lambda_2w$.
When Breaker moves, he will add $b$ edges. So,
$$ T_{i+1}\leq T_i-a\lambda_2w+\left((1+\lambda_1)^b-1\right)w.$$
To ensure Property (\ref{it:prop2}), we need to have
\begin{equation}
(1+\lambda_1)^b\leq 1+a\lambda_2 . \label{eq:cond2}
\end{equation}
In order to ensure Property (\ref{it:prop3}), we require
$$ T_0=n(1+\lambda_1)^{-bn/(a+b)-k}(1-\lambda_2)^{-an/(a+b)+k}<1 . $$
This simplifies to
\begin{equation}
1+\lambda_1>n^{\frac{a+b}{bn+k(a+b)}}
(1-\lambda_2)^{-\frac{an-k(a+b)}{bn+k(a+b)}} . \label{eq:cond3}
\end{equation}
To satisfy (\ref{eq:cond2}) and (\ref{eq:cond3}), we need
$$ n^{\frac{a+b}{bn+k(a+b)}} (1-\lambda_2)^{-\frac{an-k(a+b)}{bn+k(a+b)}}
 <1+\lambda_1\leq(1+a\lambda_2)^{1/b} . $$
Hence, we merely need to verify the existence of a $\lambda_2>0$
that gives

\begin{equation}
n^{b(a+b)}<(1+a\lambda_2)^{bn+k(a+b)} (1-\lambda_2)^{abn-kb(a+b)}.
\label{csillag}
\end{equation}

Let $\alpha$ be the unique negative root of the equation $1+x
=\exp\{x-x^2\}$. Note that $1+x\geq\exp\{x-x^2\}$ for all $x \geq
\alpha$, and $\alpha \approx -0.684$. Observe, if $(a+b)/a^2 \geq
\alpha^2 n / \ln n$, then
$$ \frac{an}{a+b}\leq\frac{\ln n}{\alpha^2a}<\frac{3\ln n}{a}\leq \frac{6ab\sqrt{n\ln n}}{(2\sqrt{ab})^{3/2}}\leq\frac{6ab\sqrt{n\ln n}}{(a+b)^{3/2}} , $$
and the statement of the lemma would be  vacuous since it merely guarantees that Maker's graph has minimum degree at least $0$.  Otherwise one may substitute
$\lambda_2=\sqrt{\frac{(a+b)\ln n}{a(a+1)n}} < -\alpha$, and use the
lower bound on $1+x$. So it is enough to see that (\ref{csillag})
holds as long as

\begin{eqnarray}
n^{b(a+b)} & < & \exp\left\{\lambda_2k(a+b)^2
+\lambda_2^2\left((b-a^2)(a+b)k-ab(a+1)n
\right)\right\} \nonumber\\
2b(a+b)\ln n & < &
k\left[(a+b)^2 \sqrt{\frac{(a+b)\ln n}{a(a+1)n}} +
\frac{(a+b)^2 \ln n}{a(a+1)n}(b-a^2)\right] \nonumber
\end{eqnarray}
\begin{eqnarray} \label{ketcsillag}
\frac{2b\sqrt{a(a+1)}}{(a+b)^{3/2}}\sqrt{n\ln n} & <  &
   k\left[ 1 +(b-a^2)\sqrt{\frac{\ln n}{a(a+1)(a+b)n}} \: \right].
\end{eqnarray}
Since it is true that
$$(b-a^2)\sqrt{\frac{\ln n}{a(a+1)(a+b)n}}\geq -\sqrt{a \ln n/n} \geq -1/2,$$
 then, by the assumption $a \leq
n/(4 \ln n)$ and the fact that $3ab \geq 2b\sqrt{a(a+1)}$,
inequality (\ref{ketcsillag}) will be satisfied if

$$\frac{6ab\sqrt{n \ln n}}{(a+b)^{3/2}} \leq k,$$
which was known.

\bizveg

\begin{lemma} \label{bdegree2}
Let $n>2a$. Breaker wins the $(a:b)$-${\rm MINDEG}(d)$ game on $K_n$ if $d>
a\left\lfloor\frac{n}{a+b}\right\rfloor$.
\end{lemma}

\noindent\textbf{Remark.} Again, this is the $(a:b)$-${\rm MINDEG}(d)$ game and, where the parameter is understood, refer to the winning strategy as \textbf{playing as MINDEG-Breaker}. \\

\noindent {\bf Proof of Lemma~\ref{bdegree2}.} In the first round
Breaker chooses a vertex, say $v$, which Maker has not touched and
chooses all of his edges to be incident to that vertex in every
round.  At the end of the game, Maker has chosen at most
$a\left\lfloor\frac{n-1}{a+b}\right\rfloor$ edges incident to  $v$.
\bizveg

\subsection{Expansion game} In the expansion game, Maker wins
if he manages to achieve that in his graph  for every pair of
disjoint sets $R$ and $S$, where $|R|=r$ and $|S|=s$, he has an edge
between $R$ and $S$. We may assume that $s \geq r$.

\begin{lemma} Maker wins the $(a:b)$-expansion game on $K_n$ with parameters
$r\leq s$ if one of the following holds:
\renewcommand{\theenumi}{\alph{enumi}}
\begin{enumerate}
   \item $2b\ln n< r\ln(a+1)$,  \label{case1}
   \item $b\ln n<r\ln(a+1)\le 2b\ln n$ and
   $s>\frac{rb\ln n}{r\ln(a+1)-b\ln n}$, \label{case2}
   \item   $n-s<\frac{nr\ln(a+1)}{b\ln n+r\ln(a+1)}$. \label{case3}
\end{enumerate}  \label{expand}
\end{lemma}

\noindent\textbf{Remark.} In referencing this game, we call it the $(a:b)$-${\rm EXP}(r,s)$ game and, where the parameters are understood, refer to the winning strategy as \textbf{playing as EXP-Maker}. \\

\noindent {\bf Proof of Lemma~\ref{expand}.} Construct a hypergraph
that has vertex set $E(K_n)$ and each hyperedge consists of the set
of edges between disjoint sets $R,S\subseteq V(K_n)$ where $|R|=r$
and $|S|=s$.  Maker attempts to occupy at least one vertex in every
hyperedge.  If he manages then he wins the expansion game.  We use
Theorem~\ref{ESB} in order to verify that each of the conditions
implies Maker's win.  Recall that the corresponding hypergraph has
$\binom{n}{r}\binom{n-r}{s}$ hyperedges, each of size $rs$.

First we prove that condition (\ref{case1}) is sufficient; i.e.,
\begin{eqnarray}
   \binom{n}{r}\binom{n-r}{s}(1+a)^{-rs/b}
   & \leq & \exp\left[(s+r)\ln n-\frac{rs}{b}\ln(a+1)\right]
   \label{esexp} \\
   & \le    & \exp\left[s\left(2\ln n-\frac{r}{b}\ln(a+1)
                        \right)\right]<1 . \nonumber
\end{eqnarray}
Assuming (\ref{case2})  we just have that
$$\binom{n}{r}\binom{n-r}{s}(1+a)^{-rs/b}
    \leq  \exp\left[(s+r)\ln n-\frac{rs}{b}\ln(a+1)\right]<1.
$$
For (\ref{case3}) we have a similar computation, we just use
$\binom{n-r}{n-r-s}=\binom{n-r}{s}$:
 $$ \binom{n}{r}\binom{n-r}{n-r-s}(1+a)^{-rs/b}
   \leq
  \exp\left[(n-s)\left(\ln n+\frac{r}{b}\ln(a+1)\right)
                       -\frac{nr}{b}\ln(a+1)\right]<1.
$$~\bizveg

\section{The $2$-diameter game}

\subsection{The $b<a$ case in the $2$-diameter game.}
First we prove that Maker wins the $\D_2(a:b)$-game if $b < a< (n/(72\ln n))^{1/3}$ and $n$ is large enough. Maker's strategy is to play the degree game with $d=\left\lceil\frac{n-1}{2}\right\rceil$ on $K_n$.  By Lemma~\ref{bdegree}, he wins the game and it is easy to check that the diameter of Maker's graph is $2$, i.e., he wins the $D_2(a:b)$-game.  Indeed, if $uv$ is not a Maker's edge for some vertices $u$ and $v$, then $|N(u)\cap N(v)|\geq n-2-2
\left(n-2-\left\lceil\frac{n-1}{2}\right\rceil\right)>0$, implying
that the intersection is non-empty. That is, in  Maker's graph, the
distance between $u$ and $v$, hence the diameter of the graph, is
at most two.

\subsection{Proof of Proposition~\ref{2-(1:1)}.}
For $n\le 3$ the statement is obvious. For $n\ge 4$ regardless of
whether Maker or Breaker starts the game, Breaker chooses an edge
not incident to that chosen by Maker. Let this edge -- the one
Breaker chooses -- be $uv$. The strategy of Breaker is: if Maker
chooses an edge incident to $u$, say $uw$, then Breaker chooses
$wv$. Similarly, if Maker chooses $vw$, then Breaker chooses $wu$.
Otherwise Breaker may take an arbitrary edge. Clearly, at the end of
the game, in Maker's graph the pair $\{u,v\}$ has distance at least three.
\bizveg

\medskip

\subsection{Proof of Theorem~\ref{2-(2:b)}}

\subsubsection{Breaker wins when $b$ is large.} We prove that Breaker wins the $\D_2(2:b)$-game for
$b=(2+\epsilon)\sqrt {n/\ln n}$. Breaker plays in two phases. In
Phase I, before his first move, he picks a vertex $v$ which has no edge in Maker's graph yet.  For $r'\le (n+b-1)/b$ rounds, he occupies as many incident
edges to $v$ as possible.  Let $u_1,\ldots,u_t$ be the list of
vertices so that Maker occupied the edge $vu_i$ before Breaker makes
his $(r'+1)^{\rm st}$ move. Trivially $t\leq 2r'+2$.  At the end of Phase I, there is no unclaimed edge incident to $v$.

In Phase II, Breaker considers $n-t-1-2r'$ disjoint sets of edges:
For each vertex $x\not\in\{v,u_1,\ldots,u_t\}$ such that neither $xv$ nor any $xu_i$ is occupied by Maker after round $r'$, define $E_x:=\{xu_1,\ldots,xu_t\}$.  By
Theorem~\ref{thmchvatal} (ii), Breaker can occupy one of these sets,
say $E_x$, when $$t \le \frac{b-1}{2}\ln (n-t-1-r'),$$ which is satisfied for
$b=(2+\epsilon)\sqrt{n/\ln n}$, if $n$ is large enough.  This forces $v$ and $x$, in Maker's graph, to be at a distance of at least $3$ from each other, i.e., Breaker won the game. \\

\subsubsection{Maker wins when $b$ is small.} We set $r$, $s$ and $c$ such that the  possible value of $b$ is maximized.  The values will be
\begin{equation}
   b=\frac{n^{1/8}}{9(\ln n)^{3/8}}, \qquad c=\frac{1}{8}, \qquad
   r=\sqrt{\frac{n\ln n}{2}} \qquad\mbox{and}\qquad
   s=\frac{n^{3/4}}{\ln n} ,
\label{eq:params}
\end{equation}
although we will not substitute these values until the end of the proof.

Maker's strategy consists of two phases. The first
one, which lasts $2nr$ rounds, uses  $2nr(b+2)$ edges, and the
second deals with the rest of the $\binom{n}{2}$ edges.  In the
first phase, Maker will play four subgames, each with a different
strategy.

Denote $\deg_B^I(x)$ to be Breaker's degree  and $\deg_M^I(x)$ to be
Maker's degree at vertex $x$ after Phase I.  In general, $\deg_B(x)$
and $\deg_M(x)$ will denote Breaker's and Maker's degrees,
respectively,
in whichever round the context indicates. \\

\noindent {\bf Phase I.} There are $2nr$ rounds in this phase. Each
of the following games is played in successive rounds.  That is,
Maker plays game $i$ in round $j$ iff $i\equiv j\pmod{4}$. A vertex
 becomes {\it high} if it achieves $\deg_B(x)\geq cn/b$ before the end of
Phase I. Note that Maker's goals are monotone properties, i.e., if
the strategy of a subgame requires Maker to occupy an edge that he already
occupied in an other subgame, then he is free to use his
edges in any way for this turn. The goals of Maker in the four games played in Phase I are the following:
\begin{itemize}
   \item {\bf Game 1. Ratio game.} If vertex $x$
    becomes high, then after this change the following relation
     will hold during the rest of Phase I: $\frac{\deg_B(x)}{\deg_M(x)}<3b$.
   \item {\bf Game 2. Degree game.}  For all vertices $x$,
   $\deg^I_M(x)\geq r$.
   \item {\bf Game 3. Expansion game.} For every pair of disjoint sets with
    $|R|\geq r$ and
   $|S|\geq s$  there is a Maker's edge between $R$
   and $S$ at the end of Phase I.
   \item {\bf Game 4. Connecting high vertices.}
   In this subgame, which lasts in the entire game not only in Phase I,
   the aim of Maker is to connect each pair of high vertices with a path of
   length at most two.
\end{itemize}

\noindent {\bf Phase II.} In the odd rounds of this phase, Maker
will connect with a path of length at most two each pair of vertices whose distance in Maker's graph is at least $3$.   As for the even turns of this
game, half of them are already dedicated to continue  Game 4, the
other half are arbitrary moves by Maker. Because Game 4 played in
the entire game, it is easier to analyze the connection of pairs of
vertices by performing it only in odd rounds.

By Game 2, after Phase I is finished, $\deg_M^I(u)\ge r$ for every vertex $u$.
By Game 3, after Phase I is finished, in Maker's graph there is an edge from the neighborhood
of $u$ into every $s$ set of
vertices, hence to all but $s$ vertices there is a path of length at
most $2$ from $u$ at the end of Phase I. The aim of Maker in Phase
II to connect the remaining pairs of vertices with a path of length
at most $2$. This is handled in Game 4 for pairs $(u,w)$ when both are
high. So in the odd rounds in Phase II we only need to connect $u$
and $w$ where either $u$ or $w$ is not a high vertex. \\

\noindent {\bf Game 1 verification.}  We can view Game 1 as a
$(2:4b)$-game.  Here, Maker plays the $(2:4b)$-degree game, for which he has
a winning strategy, provided
$d<n/(1+2b)-\frac{48b}{(2+4b)^{3/2}}\sqrt{n\ln n}$ using
Lemma~\ref{bdegree}.  Suppose, that for a vertex $x$ after becoming
high at some point in Phase I, $\deg_B(x)\geq cn/b$ but
$\deg_M(x)<\deg_B(x)/(3b)$. We claim that it is not possible, because
starting from this situation Breaker has a simple strategy to win the
min-degree game: claiming only edges incident to $x$. This way Breaker
occupies at least a $\frac{4b}{4b+2}$-fraction of the remaining edges
and his degree at the end is at least:

\begin{eqnarray}
   \lefteqn{\deg_B(x)
            +\left[n-1-\left(\deg_B(x)+\deg_M(x)\right)-(4b-1)\right]
            \frac{4b}{4b+2}} \nonumber \\
   & \geq & \deg_B(x)
            +\left[n-4b-\deg_B(x)\left(1+\frac{1}{3b}\right)
             \right]\frac{4b}{4b+2} \nonumber \\
   & \geq & (n-4b)\left(\frac{2b}{2b+1}\right)
            +\deg_B(x)\left(\frac{1}{6b+3}\right) \nonumber \\
   & \geq & n-4b-\frac n{2b+1}+\frac{8b^2}{2b+1}
            +\left(\frac{cn}{b}\right)
             \left(\frac{1}{6b+3}\right). \label{breakerdeg}
\end{eqnarray}
As long as
\begin{equation}
c\geq 36b^{3/2}\sqrt{\frac{\ln n}{n}} \label{cond1}
\end{equation}
and $n$ is large enough, we have the following:
\begin{eqnarray*}
   \frac{cn}{b(6b+3)} & \geq & \frac{36b^{3/2}}{b(6b+3)}\sqrt{n\ln n}
   > \frac{48b}{(4b+2)^{3/2}}\sqrt{n\ln n} +4b. \\
\end{eqnarray*}
So Maker occupies at most
$d=\frac{n}{2b+1}-\frac{48b}{(4b+2)^{3/2}}\sqrt{n\ln n}$ edges incident
to $x$, contradicting that he plays according to a winning strategy in
the degree game. Hence Maker wins Game 1.

\smallskip

\noindent {\bf Game 2 verification.}  Here, Maker will use a straightforward greedy strategy for the $(2:4b)$-degree game until the end of Phase
I. By the end of Phase I, Maker will have claimed $nr$  edges during
Game 2.

Maker's strategy is that in his next move, he tries to put an edge incident to a vertex whose Maker degree is the smallest.  Game 2 is played in $nr/2$ rounds, so Maker has $nr$ edges to use.  Therefore, the only reason that Maker cannot increase $\deg_M^I(x)$ (Maker's degree of $x$ in Phase I) to $r$ by the end of Phase I is that Breaker occupied more than $n-1-r$ edges incident to $x$.  More precisely, if Maker takes less than $r$ incident edges to $x$ while Breaker takes the rest, then in this case, $\deg_M^I(x)+\deg_B^I(x)=n-1$ which gives that
$\deg_B^I(x)\geq n-r$ and consequently
$\deg_B^I(x)\geq 3br$ since

\begin{equation}
n\geq 4br. \label{cond2}
\end{equation}

But then Game 1 ensures that $\deg^I_M(x)>r$, i.e. Maker wins Game 2, given that we already know that Maker wins Game 1.

\bigskip

\noindent {\bf Game 3 verification.}  Here Maker plays a {\em
virtual} $\left(2:\binom{n}{2}/(2nr)-2\right)$-expansion game
on sets of size $r$ and $s$. That is, for each $4b$ edges that
Breaker chooses, Maker also assumes that Breaker also adds a set of
$\binom{n}{2}/(2nr)-2-4b$ edges arbitrarily.  Note that we must have
that
\begin{equation}
   \binom{n}{2}/(2nr)-2\geq 4b \label{cond3c}.
\end{equation}
In this virtual game, we  apply results about  completed games,
because Phase I is finished in $2nr$ rounds. This ensures not only a
Maker's edge between each $(R,S)$ pair, but ensures that it has
occurred in Phase I.

By Lemma~\ref{expand}(a), it is sufficient
to have:
\begin{equation}
   s\geq r\geq 3 \label{cond3a}
\end{equation}
and
\begin{equation}
2\left(\frac{n}{4r}\right)\ln n < r\ln (2+1) \ \ \
\Longleftrightarrow \ \ \
   \ln n<\frac{2r^2\ln 3}{n}. \label{cond3b}
\end{equation}
(We use $n/(4r)$ in place of ``$b$'' in the lemma.)
\bigskip

\noindent {\bf Game 4 verification.}   We will see below that, when a
high vertex emerges, then from it  to   every high vertex that
emerged previously there are many paths of length $2$ consisting of
unoccupied edges.  Since Breaker can begin working immediately on
these paths, Maker must begin playing this game during Phase I.
Since high vertices can emerge at the end of Phase I, this game must
continue into Phase II.

Denote the high vertices  $x_1,x_2,\ldots,x_\ell$, indexed in the
order in which they appear. Since  Breaker occupies $2nrb$ edges, so
$4nrb$ endpoints, by the end of Phase I, hence
$$ 4nrb\geq\ell\left(\frac{cn}{b}\right), $$
implying $\ell\leq\frac{4rb^2}{c}$.  Recall that
\begin{equation}
   c=\frac{1}{8},\qquad\mbox{implying}\qquad\ell\leq 32rb^2.
   \label{eq:c18}
\end{equation}

Game 1 ensures that, for a high vertex, $x_t$, after it became high,
during Phase I $\frac{\deg_B(x_t)-4b}{\deg_M(x_t)}<3b$ holds (we
subtract $4b$ from $\deg_B(x_t)$ because the effects of Game 1 may
be delayed by a round). Hence
$\deg_B(x_t)<\frac{3b(n-1)+4b}{3b+1}<\frac{3bn}{3b+1}+\frac{1}{3}$.
Moreover, in the set of four rounds when $x_t$ becomes high,
$\deg_B(x_t)\leq\frac{cn}{b}+4b=\frac{n}{8b}+4b$.

For $j<t$ and after Round $4i$, let $Y_i(j,t)$ denote the number of
paths of length $2$ between $x_j$ and $x_t$ in which neither edge is
occupied by Breaker.   Suppose $x_t$ becomes high in round $i^*$,
then the number of such paths, available to be taken by Maker, is
$$ Y_{i^*}(j,t)\geq n-1-\deg_B(x_j)-\deg_B(x_t)
   >n\left(\frac{1}{3b+1}-\frac{c}{b}\right)-\frac{1}{3}-1-4b
  . $$
If $b$ grows slowly ($b=o(\sqrt{n})$ is sufficient), then, when
$x_t$ becomes high,
$$Y_{i^*}(j,t)\geq\frac{n}{b}\left(\frac{1}{4}-c\right)=\frac{n}{8b}.$$

We will again use a weight function argument. (See similar ideas in
Pluh\'ar \cite{PA3}.) The hypergraph is ${\mathcal H}=(V({\mathcal
H}), E({\mathcal H}))$, where $V({\mathcal H})$ is the set of paths
of length two between high vertices.\footnote{In order to avoid
confusion, we refer to a vertex or edge of ${\mathcal H}$ as {\em
hypervertex} or {\em hyperedge}, respectively.} A hyperedge (i.e.,
a winning set) is of the form $A(j,t)$, $j<t$ which consists of the
paths of length two between $x_j$ and $x_t$ which did not contain a
Breaker edge at the time $x_t$ becomes high. In each turn Maker uses his two edges to (fully) occupy paths of length $2$. Breaker, on the other hand, chooses $4b$ edges
between rounds of Game 4. If Breaker chooses edge $x_ty$, then there
are at most $\ell$ paths that he occupies to connect $x_t$ to other
high vertices. Similarly for $y$ if it also happens to be a high
vertex. As a result, this can be viewed as a $(1:8b\ell)$-game
between Maker and Breaker.

A complication with this game is that the board is {\em dynamic}. Maker
must connect the high vertices as they emerge before he knows which
vertices will become high later.  When a new vertex becomes high,
the vertex set of our hypergraph increases, as does the set of
hyperedges.  We, however, modify the weight argument to ensure that
this is not a problem.

We assign to the
hyperedge $A(j,t)$ at round $4i$ the weight
$w_i\left(A(j,t)\right)=(1+\lambda)^{-Y_i(j,t)}$ if Maker has not
connected $x_j$ and $x_t$ by a path of length $2$ and both $x_j$ and $x_t$ are
already high.  Otherwise let $w_i\left(A(j,t)\right)=0$.  The value of $\lambda >0$
will be specified later.
At any round
$4i$, $$ T_i=\sum_{j,t}w_i\left(A(j,t)\right).$$ We want to ensure
the four properties itemized below.  We take \textit{$i^{\it th}$ step},
to be the step that begins with Maker's $(4i)^{\rm th}$ move and before
Maker's $(4i+4)^{\rm th}$ move, when he, again, plays Game 4.
So, it encompasses 4 rounds of moves by both Maker and Breaker.
\renewcommand{\theenumi}{\roman{enumi}}
\begin{enumerate}
\item If Breaker wins in the $i^{\rm th}$ step, then $T_i\geq 1$,
\label{it:pr1}
\item $T_{i+1}\leq T_i+\ell(1+\lambda)^{-n/(8b)+4b}$, where $\ell$ is
the number of vertices that become high during Phase I, and if no vertices become high in the $i^{\rm th}$ step, then $T_{i+1}\leq T_i$. \label{it:pr2}
\item $T_0=0$, \label{it:pr3}
\item $\ell^2(1+\lambda)^{-n/(8b)+4b}<1$. \label{it:pr4}
\end{enumerate}
If all of the above conditions hold, then (\ref{it:pr2}) ensures $T_i$ will only increase in steps for which a vertex becomes high.  Since there are at most $\ell$ such steps, condition (\ref{it:pr4}) ensures that $T_i<1$ for all $i$.  By condition (\ref{it:pr1}), however, Breaker does not win the game.  Thus, the above is a winning strategy for Maker for Game 4.

Property (\ref{it:pr1})  trivially holds. If
\begin{equation}
   4b<cn/b , \label{cond4a}
\end{equation}
then no vertex becomes high before Maker moves in round 4 and
property (\ref{it:pr3}) is true also.

Maker follows the greedy strategy, that is he always chooses a
hypervertex of ${\mathcal H}$ (i.e., a $2$-path in $K_n$) with the
maximum weight.
Let $w$ be the weight of the largest weighted hypervertex of ${\mathcal H}$
\textbf{before} Maker makes his move. This means that Maker will reduce
the value of $T_i$ by at least $w$. When Breaker moves, he could ruin
many $2$-paths.  In particular, for any pair of high vertices, $x_j$ and
$x_t$, the weight of the hyperedge $A(j,t)$ changes by a
multiplicative factor of at most
$$ (1+\lambda)^{d_i(x_j)+d_i(x_t)}-1 , $$
where $d_i(x_t)$ is the number of edges that Breaker takes incident
to $x_t$ in the $i^{\rm th}$ step.

So, at the end of step $i$, we analyze the change in the weight.
Maker occupies a hypervertex with the maximum weight $w$ and
Breaker's moves will add weight to hyperedges, each with weight at
most $w$. Finally, we add a term for the vertices that become high
during this step. Note that if a vertex becomes high early in the
step but then the weight is increased by further play, then the
increase is still accounted for by the summation term.  Using the
following inequality, which holds for all $\alpha,\beta \geq 0$:
$$
(1+\alpha)^\beta-1\leq\exp(\alpha\beta)-1\leq\alpha\beta\exp(\alpha\beta) ,
$$
we obtain that if no vertices become high in the $i^{\rm th}$ step, then $$ T_{i+1}\leq T_i-w+\sum_{j,t}\left[(1+\lambda)^{d_i(x_j)+d_i(x_t)}-1\right]w . $$
Otherwise,
\begin{eqnarray*}
   T_{i+1} & \leq &
   T_i-w+\sum_{j,t}\left[(1+\lambda)^{d_i(x_j)+d_i(x_t)}-1\right]w
   +\ell(1+\lambda)^{-n/(8b)+4b} \\
   & \leq & T_i+\left[-1+ \sum_{j,t}\lambda\left(d_i(x_j)+d_i(x_t)\right)
   \exp\left[\lambda\left(d_i(x_j)+d_i(x_t)\right)\right]\right]w
   +\ell(1+\lambda)^{-n/(8b)+4b} \\
   & \leq & T_i+\left[-1+ \sum_{j,t}\lambda\left(d_i(x_j)+d_i(x_t)\right)
   \exp\left[\lambda(4b+1)\right]\right]w+\ell(1+\lambda)^{-n/(8b)+4b}.
\end{eqnarray*}

Since an edge is counted at most twice in the degree summation, we have $\sum_{j, t}(d_i(x_j) + d_i(x_t)) \leq 2\ell\sum_jd_i(x_j)\leq 16b\ell$, and
$$ T_{i+1}\leq T_i+\left[16\ell b\lambda\exp\left\{\lambda(4b+1)\right\}-1\right]w
+\ell(1+\lambda)^{-n/(8b)+4b} . $$

In order for condition (\ref{it:pr2}) to hold, it suffices to show
that $16b\ell\lambda\exp\left\{\lambda(4b+1)\right\}\leq 1$.  We
choose $\lambda=\frac{1}{16 b\ell}-\frac{4b+1}{(16b\ell)^2}$.
\begin{eqnarray*}
   16b\ell\lambda\exp\left\{\lambda(4b+1)\right\} & = &
   \left(1-\frac{4b+1}{16b\ell}\right)
   \exp\left\{\frac{4b+1}{16b\ell}
   -\left(\frac{4b+1}{16b\ell}\right)^2\right\} \\
   & \leq & \exp\left\{-\frac{4b+1}{16b\ell}+\frac{4b+1}{16b\ell}
   -\left(\frac{4b+1}{16b\ell}\right)^2\right\}<1 ,
\end{eqnarray*}
and so, condition (\ref{it:pr2}) is satisfied.

Finally, to find that condition (\ref{it:pr4}) holds, it is
sufficient to show that, with
$\lambda=\frac{1}{16b\ell}-\frac{4b+1}{(16b\ell)^2}$,
\begin{equation}\label{eq:lsquared}
\ell^2(1+\lambda)^{-n/(8b)+4b}\leq\exp\left\{2\ln\ell-\lambda\frac{n}{8 b}+4\lambda b\right\}<1
. \end{equation}

Since $b \geq 1$ and $\ell\geq 2$ otherwise Game 4 is irrelevant,
$$ \lambda=\frac{1}{16b\ell}-\frac{4b+1}{(16b\ell)^2}\geq
\frac{1}{16b\ell}\left(1-\frac{1}{4\ell}-\frac{1}{16b\ell}\right)
\geq\frac{1}{16b\ell}\left(\frac{27}{32}\right)>\frac{3}{64 b\ell} .
$$

Hence, because $2\leq\ell\leq 32rb^2$,
\begin{eqnarray*}
   \exp\left\{2\ln\ell-\lambda\frac{n}{8 b}+\lambda 4b\right\} & \leq &
   \exp\left\{2\ln\ell-\frac{3n}{512 b^2\ell}+\frac{3}{16\ell}\right\} \\
   & \leq & \exp\left\{2\ln(32rb^2)-\frac{3n}{16384rb^4}+\frac{3}{32}\right\}
\end{eqnarray*}
and so, in order for (\ref{eq:lsquared}) to be satisfied, it would be sufficient for
\begin{equation}
   \ln(32rb^2)-\frac{3n}{32768rb^4}+\frac{3}{64} < 0. \label{cond4b}
\end{equation}

\smallskip

\noindent {\bf Phase II verification.}
Recall that Game 4 played in Phase II also, in the even rounds.  The aim of Maker in Phase
II to connect the remaining pairs of vertices with a path of length
at most $2$. This is handled in Game 4 for pairs $(u,w)$ when both are
high. So in the odd rounds in Phase II we only need to connect $u$
and $w$ where either $u$ or $w$ is not a high vertex.

For this purpose we define the following  hypergraph:  ${\mathcal
H}=(V({\mathcal H}), E({\mathcal H}))$, where $V({\mathcal H})$ is
the set of paths of length two between vertices. For each pair of
vertices $u,w$, satisfying that at least one is not high and after Phase I
did not have a Maker's $2$-path between them, we define a hyperedge
(i.e., a winning set) $A(u,w)$, which consists of the paths of
length two between $u$ and $w$ which did not contain a Breaker edge
at the time Phase II began.

When Breaker chooses an edge $uv$ in the graph, he may take one elements
of the hyperedges $A(u,w)$ and $A(v,w)$ for all $w$. As we have seen, the
number of these type of hyperedges not more than $2s$ at the beginning of
Phase II.
Therefore, considering the hypergraph game on $\mathcal H$, it can be
viewed as a $(1:4bs)$-game between Maker and Breaker.

We have to compute the size of the winning sets.  If $u$ is high and
$w$ is not, then when Phase II begins,
$$ \deg_B(u)<\frac{3bn}{3b+1}+\frac{1}{3} \qquad\mbox{and}\qquad
   \deg_B(w)<\frac{cn}{b}+4b . $$

Therefore, for  $n$ large enough (and recalling $c=1/8$),  there
are at least $(1/4-c)n/b=n/(8b)$ available paths of length $2$
between $u$ and $w$. There are even more paths if neither $u$ nor $w$ is high.

Using Theorem~\ref{thm:es}, we see that $|E(\mathcal{H})|\leq ns/2$,
Breaker makes $4bs$ moves on each turn and Maker
(utilizing the strategy of ESB-Breaker on $\mathcal{H}$) makes $1$.  Hence, Maker
has a winning strategy if
\begin{equation}
   \frac{ns}{2}2^{-n/(32b^2s)}<1 . \label{cond5}
\end{equation}

We compile the inequalities which need to be satisfied:
(\ref{cond1}), (\ref{cond2}), (\ref{cond3c}), (\ref{cond3a}),
(\ref{cond3b}), (\ref{eq:c18}), (\ref{cond4a}), (\ref{cond4b}) and
(\ref{cond5}). Let the parameters be as defined in (\ref{eq:params}).

Conditions (\ref{cond1}), (\ref{cond2}), (\ref{cond3c}),
(\ref{cond3a}) are trivially satisfied for these values as long as
$n$ is large enough.  Conditions (\ref{cond3b}), (\ref{cond4b}) and (\ref{cond5}) determine the appropriate
values of $b$, $r$ and $s$. \\

\noindent{\bf Proving that Maker wins the ${\D}_2(2: b)$ game.}  In Phase I, by winning Game 2, Maker achieves a minimum degree $r$ in his graph. By winning Game 3 Maker ensures that for every vertex $u$ and every set $S$ with $|S|>s$ he has an edge between $S$ and the at least $r$ Maker edges incident to $u$.  I.e., after Phase I, there are at most $s$ vertices which are not within  distance $2$ of $u$. After this thinning process, in Phase II, Maker tries to connect the remaining pair of vertices with a $2$-path.

If both $\deg_B^I(u)$ and $\deg_B^I(v)$ are small then there many possibilities for Maker to build a  $2$-path between $u$ and $v$.  If only $\deg_B^I(u)$ is small, then by Game 1, $\deg_M^I(v)$ is still large and still there are many ways to build a  $2$-path between $u$ and $v$.  These path buildings are done in Phase II.

In case both $deg_B^I(u)$ and $deg_B^I(v)$ are high then the  $2$-path  between $u$ and $v$ is built in Game 4. Thus, Maker wins the ${\D}_2(2: b)$ game for $b=\frac{n^{1/8}}{9(\ln n)^{3/8}}$. \bizveg

\section{Results on the general $d$-diameter game}

For this section, let the \textit{distance between $u$ and $v$} be the length of the shortest path between vertices $u$ and $v$ in the graph induced by Maker's edges and be denoted ${\rm dist}_M(u,v)$.  For every positive integer $i$, let the \textit{$i^{\rm th}$ ball around $v$} be set of vertices of distance at most $i$ from $v$ (including $v$ itself).  Denote it by $B_i(v)$.

\subsection{Proof of Theorem~\ref{d-(1:b)} -- Maker's strategy}
Maker will play a different game in different sets of rounds. For
$k=1,\ldots,\lceil d/2\rceil-1$, Maker will play Game $k$ in round
$i$ when $i\equiv k\pmod{\lceil d/2\rceil}$.  In rounds $i$, where
$i\equiv 0\pmod{\lceil d/2\rceil}$, Maker plays game $\lceil
d/2\rceil$. We first describe the individual subgames in detail. \\

We define several variables for the games.  There is a sequence $1=r_0,r_1,\ldots,r_{\lceil d/2\rceil-1}$ such that
\begin{eqnarray*}
   r_1 & = & \frac{n}{\lceil d/2\rceil b}
             \left(1-7\sqrt{\frac{\lceil d/2\rceil b\ln n}
                                 {n}}\right), \mbox{ and} \\
   r_i & = & \frac{nr_{i-1}\ln 2}
                  {\lceil d/2\rceil b\ln n+r_{i-1}\ln 2}
             -\sum_{j=0}^{i-1}r_j \ \ \text{ for }
   \qquad i=2,\ldots,\lceil d/2\rceil-1 .
\end{eqnarray*}
Let
$$ \beta:=\left(\frac{2n\ln 2}{\ln n}\right)^{1/\lceil
d/2\rceil}\quad\mbox{and}\quad b:=\frac{n\ln 2}{\lceil
d/2\rceil\ln n}\beta^{-1} .
$$

After defining the games, we immediately prove that Maker can win them, as the justifications are straightforward. \\

\noindent \textbf{Game 1.} This game is a $\left(1:\lceil d/2\rceil b\right)$-${\rm MINDEG}(r_1)$ game.  That is, Maker plays to achieve minimum degree at least $r_1$ in the graph of edges that he occupies.  By Lemma~\ref{bdegree}, Maker is able to ensure that the minimum degree is at least
\begin{eqnarray*}
   \lefteqn{\frac{n}{1+\lceil d/2\rceil b}
                   \left(1-\frac{6\lceil d/2\rceil b}
                                {\sqrt{1+\lceil d/2\rceil b}}
                           \sqrt{\frac{\ln n}{n}}\right)} \\
   & \geq & \frac{n}{\lceil d/2\rceil b}
            \left(1-\frac{6\sqrt{\lceil d/2\rceil b}}{\sqrt{1+\lceil d/2\rceil b}}
                   \sqrt{\frac{\lceil d/2\rceil b\ln
                   n}{n}}\right) \\
   & \geq & \frac{n}{\lceil d/2\rceil b}
            \left(1-6\sqrt{\frac{\lceil d/2\rceil b\ln n}
                                {n}}\right)=r_1 ,
\end{eqnarray*}
as long as $n$ is large enough. \\

\noindent \textbf{Game 2 to Game $\bf \lceil d/2\rceil-1$.} Game $i$, for $i=2,\ldots,\lceil d/2\rceil-1$, is a $\left(1:\lceil d/2\rceil b\right)$-${\rm EXP}(r_{i-1},n-r_i)$ game.  That is, Maker's aim is to have one of his edges between any set of size $r_{i-1}$ and any set of size $n-r_i$.  So, by Lemma~\ref{expand} (\ref{case3}), if Maker's $(i-1)^{\rm st}$ neighborhood is of size at least $r_{i-1}$, then he can ensure that the $i^{\rm th}$ neighborhood is of size at least
$$ \frac{nr_{i-1}\ln 2}{\lceil d/2\rceil b\ln n+r_{i-1}\ln 2}-\sum_{j=1}^{i-1}r_j-1 . $$ \\

\noindent \textbf{Game $\bf \lceil d/2\rceil$.} The game in this case depends on the parity of $d$.

For $d$ even, the game is a $\left(1:\lceil d/2\rceil b\right)$-${\rm EXP}\left(r_{\lceil d/2\rceil-1},\lceil n/2\rceil-1\right)$ game. In that case, the $(d/2)^{\rm th}$ neighborhood of every vertex, in Maker's graph, is at least $\lfloor n/2\rfloor+1$, which implies the diameter of Maker's graph is at most $d$.

For $d$ odd, the game is a
$\left(1:\lceil d/2\rceil b\right)$-${\rm EXP}(r_{\lceil d/2\rceil-1},n-r_{\lceil d/2\rceil-1})$ game.  This ensures that Maker has an edge between the $(d-1)/2$ neighborhoods of each pair of vertices.

To prove that Maker wins this game, in both cases, we apply Lemma~\ref{expand} (\ref{case1}).  In checking the conditions of Lemma~\ref{expand} (\ref{case1}), the first is the trivial condition that $r\leq s$.  In either case, if $r>s$ then the $\left(\lceil d/2\rceil-1\right)^{\rm st}$ neighborhoods of each vertex is at least $\lfloor n/2\rfloor+1$, giving that the diameter is less than $d$ already.

For the nontrivial condition of Lemma~\ref{expand} (\ref{case1}),
\begin{equation} 2\lceil d/2\rceil b\ln n<r_{\lceil d/2\rceil-1}\ln 2. \label{eq:nontriv} \end{equation}
To prove that Maker wins this game we must justify (\ref{eq:nontriv}).

\noindent\textbf{Checking the conditions imposed by the subgames.} The following claim merely proves a statement about the constants we defined at the beginning of the section.

\begin{claim}
If $n$ is large enough, then for $i=1,\ldots,\lceil d/2\rceil-1$ we have
   \begin{equation}\label{claimeq} \left(1-6\beta^{-1/2}\right)^i
      \frac{\ln n}{\ln 2}\beta^i\leq r_i\leq\frac{\ln n}{\ln 2}\beta^i . \end{equation}
   \label{claimri}
\end{claim}

\noindent{\textbf{Proof of Claim \ref{claimri}.}} We will prove this
by induction on $i$. For $i=1$  using the definition of $r_1$ we easily see
that \eqref{claimeq} holds:
\begin{eqnarray*}
   r_1  =  \frac{n}{\lceil d/2\rceil b}
             \left(1-6\sqrt{\frac{\lceil d/2\rceil b\ln n}
                                 {n}}\right)
    = \frac{\ln n}{\ln 2}\beta\left(1-6\sqrt{\frac{\ln 2}{\beta}}\right) .
\end{eqnarray*}
Let us assume the bounds hold for $r_{i-1}$ whenever
$i\in\{2,\ldots,\lceil d/2\rceil-1\}$.  Now we compute $r_i$:
\begin{eqnarray*}
   r_i = \frac{nr_{i-1}\ln 2}
                  {\lceil d/2\rceil b\ln n+r_{i-1}\ln 2}
             -\sum_{j=0}^{i-1}r_j
    =  \frac{nr_{i-1}\beta}{n+r_{i-1}\beta}
             -\sum_{j=1}^{i-1}r_j-1 \leq r_{i-1}\beta .
\end{eqnarray*}

So by the inductive
hypothesis, $r_i\leq\frac{\ln n}{\ln 2}\beta^i$.  As to the lower
bound, we note that $n$ large enough ensures that $r_{i-1}$ is the largest among $r_1,\ldots,r_{i-1}$.
\begin{eqnarray*}
   r_i\geq r_{i-1}\beta\left[\frac{n}{r_{i-1}\beta+n}\right]-ir_{i-1}\geq r_{i-1}\beta\left[1-\frac{r_{i-1}\beta}{n}-i\beta^{-1}\right] .
\end{eqnarray*}
To bound this,
$$ i\beta^{-1}\leq (d/2)\left(\frac{\ln n}{2 n\ln 2}\right)^{1/d}
   \beta^{-1/2}\leq 3\beta^{-1/2} , $$
as long as $d\leq c\ln n/\ln\ln n$ for some constant $c$ and $n$ sufficiently large. Also, the bound $r_{i-1}\leq\frac{\ln n}{\ln 2}\beta^{i-1}$ gives
$$ \frac{r_{i-1}\beta}{n}\leq\frac{\ln n}{n\ln 2}\beta^i
   \leq\frac{\ln n}{n\ln 2}\beta^{\lceil d/2\rceil-1}=2\beta{-1}
   \leq 3\beta^{-1/2} , $$
for any $d\geq 1$.   So, we combine these and see that $r_i\geq
r_{i-1}\beta\left(1-6\beta^{-1/2}\right)$. By the inductive
hypothesis, $r_i\geq \left(1-6\beta^{-1/2}\right)^i\frac{\ln n}{\ln
2}\beta^i$ and the proof is finished. \bizveg \\

Claim 1 proves that all of the $\lceil d/2\rceil$ subgames are won by Maker and so Maker can construct a path of length at most $d$ between any pair of vertices. In particular, our values of $b$ and $r_{\lceil d/2\rceil-1}$ ensure that inequality (\ref{eq:nontriv}) is satisfied. Indeed,
  for any pair of vertices $u,v$ Maker can achieve that the $\lfloor d/2\rfloor$-neighborhoods of both vertices have sizes at least $n/2$, therefore they are intersecting and their distance in Maker's graph will be at most $d$.

This completes the proof that Maker's strategy  is successful and
$$ b=\frac{n\ln 2}{\lceil d/2\rceil \ln n}\beta^{-1}>\frac{1}{2d}\left(\frac{n}{\ln n}\right)^{1-1/\lceil d/2\rceil} . $$ \\

\subsection{Proof of Theorem~\ref{d-(1:b)} -- Breaker's strategy, $a=1$}

We will show that Breaker has a winning strategy for $a=1$, $b \geq
4d^{1/(d-1)}n^{1-1/(d-1)}$ and $n$ is sufficiently large.

In Breaker's first move, he will choose an edge between $u$ and $v$.
(If Maker goes first, Breaker will make sure that neither $u$ nor
$v$ is incident to Maker's first edge.)  We give a strategy for Breaker to ensure that ${\rm dist}_M(u,v)>d$.  Then at each move Breaker will play two
roles. He will use $d^{1/(d-1)}n^{1-1/(d-1)}$ edges to ensure the
maximum Maker degree is small and then will use the remaining edges
to ensure that every vertex in the $i^{\rm th}$ Maker's neighborhood
of $u$ is adjacent by a Breaker edge to every vertex in the $j^{\rm
th}$ Maker's neighborhood of $v$ for $j=0,\ldots,d-1-i$.  If Breaker
succeeds in his plan, then  since Maker can have no edge between the
$i^{\rm th}$ ball around $u$ and the $(d-i-1)^{\rm st}$ ball around
$v$, it is impossible for Maker to occupy any path of length at most
$d$ between $u$ and $v$, and Breaker wins.

Breaker will play as MINDEG-Maker with $a_{{\rm MINDEG}}=b_1:=d^{1/(d-1)}n^{1-1/(d-1)}$ and $b_{{\rm MINDEG}}=1$.  By Lemma~\ref{bdegree}, Breaker (i.e., MINDEG-Maker) has a
winning strategy to keep the maximum degree in Maker's graph,
$\Delta$ (i.e., minimum degree in Breaker's graph at least $n-1-\Delta$), at most
\begin{eqnarray}
\Delta & \leq & n-1-\left(\frac{b_1n}{b_1+1}
+\frac{6b_1}{(b_1+1)^{3/2}}\sqrt{n\ln n}\right) \nonumber \\
& \leq & \frac{n}{b_1} +6\sqrt{\frac{n\ln n}{b_1}} \label{delb} \\
& \leq & \left(\frac{n}{d}\right)^{1/(d-1)}
+6d^{-1/(2d-2)}\sqrt{n^{1/(d-1)}\ln n} \nonumber \\
& \leq & \left(\frac{n}{d}\right)^{1/(d-1)} \exp\left\{6d^{1/(2d-2)}
\sqrt{\frac{\ln n}{n^{1/(d-1)}}}\right\} . \nonumber
\end{eqnarray}

The remaining number of edges used by Breaker in a turn  is
$b_2:=b-b_1$. We will guarantee that, for some $k\in\{0,\ldots,d\}$,
the intersection of the balls around $u$ and $v$, namely $B_k(u)\cap
B_{d-k}(v)$, is empty. Let $xy$ be the edge chosen by Maker in the
preceding move.

For any pair of disjoint vertex sets $X$ and $Y$, $E(X,Y)$ denotes the set of edges with one endpoint in $X$ and the other in $Y$.  There are 2 cases in Breaker's strategy.
\begin{enumerate}
   \item \textbf{Vertex $x$ is as close to $u$ as $y$ is and vertex $y$ is as close to $v$ as $x$ is.}\\
   \indent Let $i=\dist(x,u)$ and $j=\dist(y,v)$.  If $i\leq d-1$, then Breaker will choose each edge in
   $\bigcup_{k=0}^{d-i-1}E\left(N_k(x),B_{d-i-1-k}(v)\right)$.  If $j\leq d-1$, then Breaker will choose each edge in
   $\bigcup_{k=0}^{d-j-1}E\left(N_k(y),B_{d-j-1-k}(u)\right)$. \label{breakeritem1} \\

   \item \textbf{Vertex $x$ is closer than $y$ to both $u$ and $v$.}\\
   \indent Let $i=\dist(x,u)$ and $j=\dist(x,v)$.
   If $i\leq d-2$, then Breaker will choose each edge in
   $\bigcup_{k=0}^{d-i-2}E\left(N_k(y),B_{d-i-2-k}(v)\right)$.  If $j\leq d-2$, then Breaker will choose each edge in
   $\bigcup_{k=0}^{d-j-2}E\left(N_k(y),B_{d-j-2-k}(u)\right)$. \label{breakeritem2}

\end{enumerate}~\\

Suppose for a  contradiction, that Maker chooses the edge $\{x',y'\}$ which results in a path of length at most $d$ from $u$ to $v$.  Let $x'$ be closer to $u$ than $y'$ and $x'\in N_i(u)$ and $y'\in N_j(v)$.  By the supposition, $i+j\leq d-1$.  Before Maker chooses $\{x',y'\}$, we look for the first point in the game at which $x'$ is distance $i$ from $u$ and $y'$ is distance $j$ from $v$.

By symmetry we may assume, without loss of generality, that $x'$ becoming distance $i$ from $u$ does not happen before $y'$ becomes distance $j$ from $v$.  When $x'$ becomes a member of $N_i(u)$, it means that there is another edge $\{x'',y''\}$ that Maker chose in a path of length $i$ from $u$ to $x'$.  Without loss of generality, we may choose $x''$ to be closer to $u$ than $y''$, denoting $x''\in N_{\ell}(u)$.  Hence, $x'\in N_{i-\ell}(x'')$.

If $\{x'',y''\}$ is in case (\ref{breakeritem1}), then Breaker will occupy all edges in $E\left(N_k(x''),B_{d-\ell-1-k}(v)\right)$ for any $k\leq d-\ell-1$.  Since $i\leq d-1$, $i-\ell\leq d-\ell-1$ and we may choose $k=i-\ell$.  In particular, Breaker occupies every edge in
$E\left(\{x'\},B_{d-i-1}(v)\right)$, since $d-i-1=d-\ell-1-(i-\ell)$.
Since $j\leq d-i-1$, $y'\in N_j(v)\subseteq B_{d-i-1}(v)$ and $\{x',y'\}$ was chosen by Breaker, so it is not available for Maker to take later.

If $\{x'',y''\}$ is in case (\ref{breakeritem2}), then Breaker will occupy all edges in \linebreak $E\left(N_k(y''),B_{d-\ell-2-k}(v)\right)$ for any $k\leq d-\ell-2$.  Since $i\leq d-1$, $i-\ell-1\leq d-\ell-2$ and we may choose $k=i-\ell-1$.  In particular, Breaker occupies every edge in
$E\left(\{x'\},B_{d-i-1}(v)\right)$, since $d-i-1=d-\ell-2-(i-\ell-1)$.
Since $j\leq d-i-1$, $y'\in N_j(v)\subseteq B_{d-i-1}(v)$ and $\{x',y'\}$ was chosen by Breaker, so it is not available for Maker to take later. \\

Having thus proven that the above strategy by Breaker will ensure diameter at least $d+1$, we compute the value of $b_2$.

Consider the number of edges that Breaker takes in case (\ref{breakeritem1}).  The total number of edges needed is at most:
\begin{eqnarray}
\lefteqn{\sum_{k=0}^{d-i-1}\Delta^k \sum_{\ell=0}^{d-i-1-k}\Delta^{\ell}
+\sum_{k=0}^{d-j-1}\Delta^k \sum_{\ell=0}^{d-j-1-k}\Delta^{\ell}} \nonumber \\
& = &
\frac{1}{\Delta-1}\sum_{k=0}^{d-i-1}\left(\Delta^{d-i}-\Delta^k\right)
+ \frac{1}{\Delta-1}\sum_{k=0}^{d-j-1}\left(\Delta^{d-j}-\Delta^k\right)
\nonumber \\
& = & \frac{1}{(\Delta-1)^2} \left[(d-i)\Delta^{d-i+1}-(d-i+1)\Delta^{d-i}+1\right]
\label{del} \\
& & +\frac{1}{(\Delta-1)^2} \left[(d-j)\Delta^{d-j+1}-(d-j+1)\Delta^{d-j}+1\right] . \nonumber
\end{eqnarray}

Case (\ref{breakeritem2}) gives a similar expression
\begin{eqnarray}
\lefteqn{\frac{1}{(\Delta-1)^2}\left[(d-i-1)\Delta^{d-i}-(d-i)\Delta^{d-i-1}+1\right]} \label{del2} \\
& & +\frac{1}{(\Delta-1)^2}\left[
(d-j-1)\Delta^{d-j}-(d-j)\Delta^{d-j-1}+1\right] . \nonumber
\end{eqnarray}

Observe that, since Breaker's strategy is successful at preventing a length $d$ path from $u$ to $v$, $i+j\geq d$ for case (\ref{breakeritem1}) and $i+j\geq d+1$ for case (\ref{breakeritem2}). Claim~\ref{claim2} is a general statement for integers:
\begin{claim}\label{claim2}  Let $\Delta, n$ be positive integers. Then for
$k=1,\ldots,n-1$,
\begin{eqnarray*}
   f(k) & := & \left[(n-k)\Delta^{n-k+1}-(n-k+1)\Delta^{n-k}+1\right]
   +\left[k\Delta^{k+1}-(k+1)\Delta^k+1\right] \\
   & \leq & (n-1)\Delta^n-n\Delta^{n-1}+1+(\Delta-1)^2.
\end{eqnarray*}
\end{claim}

\noindent {\bf Proof of Claim \ref{claim2}.} By symmetry,
$f(k)=f(n-k)$. So it is sufficient to show that $f(k)\geq f(k+1)$
for $k\in\{1,\ldots,\lfloor (n-1)/2\rfloor\}$.
\begin{eqnarray*}
   \lefteqn{f(k)-f(k+1)} \\
   & = & \left[(n-k)\Delta^{n-k+1}-(n-k+1)\Delta^{n-k}+1\right]
   +\left[k\Delta^{k+1}-(k+1)\Delta^k+1\right] \\
   & & -\left[(n-k-1)\Delta^{n-k}-(n-k)\Delta^{n-k-1}+1\right]
   -\left[(k+1)\Delta^{k+2}-(k+2)\Delta^{k+1}+1\right] \\
   & = &
   \Delta^{n-k-1}\left[(n-k)\Delta^2-2(n-k)\Delta+(n-k)\right]
   \\
   & & -\Delta^k\left[(k+1)\Delta^2-2(k+1)\Delta+(k+1)\right] \\
   & = & (\Delta-1)^2\Delta^k(n-k)
   \left[\Delta^{n-2k-1}-\frac{k+1}{n-k}\right] .
\end{eqnarray*}

Since $k\leq\lfloor (n-1)/2\rfloor$,
$$ \Delta^{n-2k-1}\geq 1\geq\frac{k+1}{n-k} . $$
Thus, it must be the case that $f(k)$ attains its maximum at $k=1$
and $k=n-1$ and the claim is proven. \bizveg \\

Therefore, in case (\ref{breakeritem1}) (where $i+j\geq d$), the term in (\ref{del}) is bounded above by
$$ \frac{1}{(\Delta-1)^2}\left[(2d-i-j-1)\Delta^{2d-i-j}-(2d-i-j)\Delta^{2d-i-j-1}+1+(\Delta-1)^2\right] . $$
In case (\ref{breakeritem2}) (where $i+j\geq d+1$), the term in (\ref{del2}) is bounded above by
$$ \frac{1}{(\Delta-1)^2}\left[(2d-i-j-3)\Delta^{2d-i-j-2}-(2d-i-j-2)\Delta^{2d-i-j-3}+1+(\Delta-1)^2\right] . $$

It can be shown that, for $\Delta\geq 3$, both expressions are maximized by making $i+j$ as small as possible so in any case, $b_2$ is bounded above by
\begin{equation} \label{delfinal}
\frac{1}{(\Delta-1)^2}\left[(d-1)\Delta^{d}-d\Delta^{d-1}+1+(\Delta-1)^2\right] .
\end{equation}

As long as $\Delta\geq d$, inequality (\ref{delfinal}) can be bounded so that $b_2\leq d\Delta^{d-2}\exp\{1/\Delta\}$ whence
\begin{eqnarray*}
   b_2\leq d\Delta^{d-2}\exp\left\{\frac{1}{\Delta}\right\}
   & \leq & d\left(\frac{n}{d}\right)^{\frac{d-2}{d-1}}
   \exp\left\{6(d-2)d^{1/(2d-2)} \sqrt{\frac{\ln n}{n^{1/(d-1)}}}
   +\left(\frac{d}{n}\right)^{\frac{1}{d-1}}\right\} \\
   & \leq & d^{1/(d-1)}n^{1-1/(d-1)}
   \exp\left\{7d^{\frac{2d-1}{2d-2}}\sqrt{\frac{\ln n}{n^{1/(d-1)}}}
   \right\} .
\end{eqnarray*}
As long as $d\leq\frac{\ln n}{3\ln\ln n}$ and $n$ is large enough,
we have that $$
b_2\leq d^{1/(d-1)}n^{1-1/(d-1)}\exp\left\{7d^{\frac{2d-1}{2d-2}}\sqrt{\frac{\ln
n}{n^{1/(d-1)}}}\right\}\leq 2d^{1/(d-1)}n^{1-1/(d-1)} . $$

\noindent Hence, we only used $b_1+b_2\geq
3d^{1/(d-1)}n^{1-1/(d-1)}$ in order to ensure a Breaker win. \\

\subsection{Proof of Theorem~\ref{d-(1:b)} -- Breaker's strategy, $a\geq 2$}

For the sake of simplicity, we deal only with the case $a=2$; the
general case is very similar.
We will show that Breaker can win if $b \geq 4n^{1-1/d}$ and will do
so by simply playing the degree game.  As we have seen in
(\ref{delb}), Breaker can ensure that
 $$  \Delta  \leq  \frac{2n}{b} +12\sqrt{\frac{n\ln n}{b}} \\
    \leq  \frac{n}{4n^{1-1/d}}+12\sqrt{\frac{n\ln n}{4n^{1-1/d}}}
   \\
    \leq  n^{1/d}\left(\frac{1}{2}+6\sqrt{n^{-1/d}\ln n}\right) \\
    \leq  (2/3)n^{1/d} ,$$
as long as $n$ is large enough.

With $\Delta$ being the maximum degree, for any vertex $v$,
   $$|B_d(v)|  \leq  1+\sum_{i=0}^{d-1}\Delta(\Delta-1)^i
    =  1+\Delta\frac{(\Delta-1)^d-1}{(\Delta-1)-1}
    \leq
   1+\left(1+\frac{2}{\Delta-2}\right)\left(\Delta^d-1\right) <2\Delta^d
   ,$$
as long as $\Delta\geq 4$.

But then, using the upper bound for $\Delta$ and recalling $d\geq
2$,
$$
|B_d(v)|<2\Delta^d<2\left(\frac{2}{3}n^{1/d}\right)^d\leq\frac{8}{9}n
. $$ So, for any vertex $v$, there is at least one vertex of
distance greater than $d$ from it. Note that in order for Breaker to
win this degree game, we need to have that $b\leq n/(4\ln n)$.
Since $b=4n^{1-1/d}$, this holds for $d\leq \ln n/(2\ln\ln n)$ and
$n$ large enough. \bizveg

\section{Thanks}
We thank Angelika Steger for uncovering a flaw in a previous version
of the manuscript and J\'ozsef Beck for reviewing it and making
helpful comments. We are indebted to the referees for their careful
reading and numerous suggestions which improved the presentation of
the paper.

\end{document}